\documentclass[12pt]{amsart}

\usepackage{graphicx}
\usepackage{amsmath}
\usepackage{amssymb}

\usepackage{maple2e}
 
\DefineParaStyle{Maple Output}
\DefineParaStyle{Maple Output}
\DefineParaStyle{Maple Plot}
\DefineCharStyle{2D Math}
\DefineCharStyle{2D Output}

\newtheorem{thm}{Theorem}[section]
\newtheorem{lem}[thm]{Lemma}%[section]
\newtheorem{prop}[thm]{Proposition}%[section]
\newtheorem{cor}[thm]{Corollary}%[section]
\theoremstyle{definition}

\theoremstyle{remark}
\newtheorem{remark}{Remark}[section] % \renewcommand{\theremark}{}

\theoremstyle{plain}

\numberwithin{equation}{section}

\def\CC{{\mathbb C}}
\def\EE{{\mathbb E}}

\def\HH{{\mathbb H}}

\def\QQ{{\mathbb Q}}
\def\RR{{\mathbb R}}
\def\TT{{\mathbb T}}

\def\ZZ{{\mathbb Z}}

\def\e{\mathrm{e}}
\def\i{\mathrm{i}}

\def\L{\operatorname{L{}}}

\def\SL{\operatorname{SL}}

\def\PSL{\operatorname{PSL}}

\def\T{\operatorname{T{}}}

\def\meas{\operatorname{meas}}

\def\Prob{\operatorname{Prob}}

\def\GamG{\Gamma\backslash G}

\def\scrF{{\mathcal F}}

\def\Re{\operatorname{Re}}
\def\Im{\operatorname{Im}}

\begin{document}

\title{Distribution modulo one and Ratner's theorem}
\author{Jens Marklof}
\address{School of Mathematics, University of Bristol,
Bristol BS8 1TW, U.K.} 
\address{{\tt j.marklof@bristol.ac.uk}}
\date{\today}

%\begin{abstract}
%\end{abstract}

\maketitle

\tableofcontents

\section{Introduction}

Measure rigidity is a branch of ergodic theory that has recently contributed to the solution of some fundamental problems in number theory and mathematical physics. Examples are proofs of quantitative versions of the Oppenheim conjecture \cite{Eskin98}, related questions on the spacings between the values of quadratic forms \cite{Eskin05,pairI,pairII}, a proof of quantum unique ergodicity for certain classes of hyperbolic surfaces \cite{Lindenstrauss06}, and an approach to the Littlewood conjecture on the nonexistence of multiplicatively badly approximable numbers \cite{Einsiedler06}. 

In these lectures we discuss a few simple applications of one of the central results in measure rigidity: Ratner's theorem. We shall investigate the statistical properties of certain number theoretic sequences, specifically the fractional parts of $m\alpha$, $m=1,2,3,\ldots$, (a classical, well understood problem) and of $\sqrt{m\alpha}$ (as recently studied by Elkies and McMullen \cite{Elkies04}). By exploiting equidistribution results on a certain homogeneous space $\GamG$, we will show that the statistical properties of these sequences can exhibit significant deviations from those of independent random variables. The ``randomness'' of other, more generic sequences such as $m^2\alpha$ and $2^m \alpha$ mod 1 has been studied extensively. We refer the interested reader to the review \cite{energy}, and recommend the papers \cite{Rudnick98,Rudnick02} as a first read.

These notes are based on lectures presented at the Institute Henri Poincar\'e Paris, June 2005, and at the summer school `Equidistribution in number theory', CRM Montr\'eal, July 2005. The author gratefully acknowledges support by an EPSRC Advanced Research Fellowship.

\section{Randomness of point sequences mod 1}

Consider an infinite triangular array of numbers on the circle $\TT=\RR/\ZZ$ (which we represent as the unit interval $[0,1)$ with its endpoints identified),
\begin{equation}\label{triarray}
\begin{array}{ccccc}
\xi_{11} &  &  &  &\\
\xi_{21} & \xi_{22} &  & &\\
\vdots & \vdots & \ddots & & \\
\xi_{N1} & \xi_{N2} & \ldots & \xi_{NN} & \\
\vdots & \vdots & & & \ddots  
\end{array}
\end{equation}
We assume that each row is ordered, i.e., $\xi_{Nj}\leq \xi_{N(j+1)}$, and are interested in quantifying statistical properties of the $N$th row as $N\to\infty$. To simplify notation we will from now on drop the index $N$, and simply write $\xi_j$ instead of $\xi_{Nj}$.

As we shall see later, many interesting statistical properties of a sequence on $\TT$ can be derived from the knowledge of the number of elements in small subintervals of $\TT$. Let $\chi$ denote the characteristic function of the interval $[-\tfrac12,\tfrac12)\subset\RR$. That is, $\chi(x)=1$ if $-\tfrac12\leq x < \tfrac12$ and $=0$ otherwise. The characteristic function of the interval $[x_0-\tfrac{\ell}{2},x_0+\tfrac{\ell}{2})+\ZZ\subset\TT$ ($\ell\leq 1$) can be represented as
\begin{equation}
\chi_\ell(x)=\sum_{n\in\ZZ} \chi\bigg(\frac{x-x_0+n}{\ell}\bigg).
\end{equation}
The sum over $n$ makes sure $\chi_\ell$ is periodic.
The number of elements $\xi_j$ in the interval are therefore
\begin{equation} \label{SN}
S_N(\ell)=\sum_{j=1}^N \chi_\ell(\xi_j).
\end{equation}

We will always assume that the rows in our triangular array become uniformly distributed mod one.  This means that for every $x_0,\ell$,
\begin{equation}
\lim_{N\to\infty} \frac1N S_N(\ell) = \ell,
\end{equation}
i.e., the proportion of elements in any given interval is asymptotic to the interval length $\ell$.\footnote{If a sequence $\{\xi_j\}$ fails to be uniformly distributed but still has a resonable limiting density $\rho$, we may rescale the $\xi_j$ to obtain a uniformly distributed sequence. This is done as follows. Suppose for every $x_0,\ell$
$
	\lim_{N\to\infty} \frac1N S_N(\ell) = \int_{x_0-\ell/2}^{x_0+\ell/2} \rho(x) dx ,
$
where the integrated density
$
	N(x)= \int_0^x \rho(x') dx'
$
is continuous and strictly increasing. We rescale the sequence $\{\xi_j\}$ by setting 
$
	\tilde \xi_j := N(\xi_j) .
$
Note that $N(\xi_j)\in[0,1)$ for $\xi_j\in[0,1)$. The new sequence $\{\tilde\xi_j\}$ is indeed uniformly distributed modulo one (exercise).}

The aim is now to characterize the different degrees of ``randomness'' of the deterministic sequence $\{\xi_j\}$ in terms of their distribution in very small intervals with random center $x_0$. A convenient length scale is the average spacing between elements, which is $1/N$. We set
\begin{equation}
L=N\ell.
\end{equation} 
We assume $x_0$ is a random variable uniformly distributed on $\TT$ with respect to Lebesgue measure $dx_0$. We will denote expectation values by
\begin{equation}
\langle \ldots \rangle = \int_0^1 \ldots dx_0 .
\end{equation}
It is easy to work out the expectation value for the number of elements in a random interval of size $\ell$, 
\begin{equation}
\langle S_N(\ell)\rangle= L .
\end{equation}
The variance is much less trivial. Let us begin by deriving a convenient representation in terms of the pair correlation density.
We have for the mean square (the ``number variance'')
\begin{equation}
\Sigma^2_N(L):=\langle [S_N(\ell)-L]^2\rangle = \langle S_N(\ell)^2 \rangle - L^2
\end{equation}
and
\begin{equation}
\begin{split}
\langle S_N(\ell)^2 \rangle 
& = \sum_{i,j=1}^N \sum_{m,n\in\ZZ}\int_0^1  \chi\bigg(\frac{\xi_i-x_0+m}{\ell}\bigg) \chi\bigg(\frac{\xi_j-x_0+n}{\ell}\bigg) dx_0 \\
& = \sum_{i,j=1}^N \sum_{m\in\ZZ}\int_{\RR}  \chi\bigg(\frac{\xi_i-x_0+m}{\ell}\bigg) \chi\bigg(\frac{\xi_j-x_0}{\ell}\bigg) dx_0 \\
& = \ell \sum_{i,j=1}^N \sum_{m\in\ZZ} \Delta\bigg(\frac{\xi_i-\xi_j+m}{\ell}\bigg) 
\end{split}
\end{equation}
where
\begin{equation}
\Delta(x)=\int_\RR \chi(x-x_0)\chi(x_0) dx_0 =\max\{ 1-|x|, 0\}.
\end{equation}
Now the diagonal terms $i=j$ in the above double sum can be easily evaluated. We have 
\begin{equation}
\ell \sum_{i=j=1}^N \sum_{m\in\ZZ} \Delta\bigg(\frac{m}{\ell}\bigg)=\ell \Delta(0) =\ell
\end{equation}
for $\ell<1$.

The {\em pair correlation function} (also called {\em two-point correlation function}) for the sequence $\{\xi_j\}$ is defined by
\begin{equation}
R_N^2(L,\psi)=\frac1N \sum_{i\neq j=1}^N \sum_{m\in\ZZ} \psi\bigg(\frac{\xi_i-\xi_j+m}{\ell}\bigg) ,
\end{equation} 
where $\psi$ is taken from a class of sufficiently nice test functions (e.g. continuous with compact support such as $\Delta$).
With the above calculation we therefore have the identity
\begin{equation}
\Sigma^2_N(L)=L-L^2+ L R_N^2(L,\Delta) .
\end{equation}
This says that the asymptotic analysis of the pair correlation density will give us information on the number variance.

Note that by the Poisson summation formula
\begin{equation}
\sum_{m\in\ZZ} f(m) = \sum_{n\in\ZZ} \widehat f(n), 
\end{equation}
where
\begin{equation}
\widehat f(y) = \int_\RR f(x) e(xy) dy, \qquad e(x):=\exp 2\pi\i x,
\end{equation}
we have
\begin{equation} \label{R22}
R_N^2(L,\psi)=\frac{L}{N^2} \sum_{i\neq j=1}^N \sum_{n\in\ZZ} \widehat\psi\bigg(\frac{L n}{N}\bigg) e\big(n(\xi_i-\xi_j)\big) .
\end{equation} 
Here $\psi$ can be any function with absolutely convergent Fourier series (e.g. $\Delta$).

\subsection{Distribution of gaps}

A popular statistical measure is the distribution of gaps 
\begin{equation}
s_j=N(\xi_{j+1}-\xi_j) \qquad (j=1,\ldots, N, \; \xi_{N+1}:=\xi_1+1)
\end{equation}
between consecutive elements (recall the $\xi_j$ form an ordered sequence on $\TT$).
We have multiplied the actual gap $\xi_{j+1}-\xi_j$ by $N$, which means we are measuring spacings in units of the average gap $1/N$.

The gap distribution of the sequence $\xi_1,\ldots,\xi_N$ is defined as
\begin{equation}
	P_N(s)=\frac1N \sum_{j=1}^N \delta(s-s_j)
\end{equation}
where $\delta$ is a Dirac mass at the origin. The question we will investigate is whether $P_N(s)$ has a limiting distribution $P(s)$. That is, does there exist a probability density $P(s)$ such that for every bounded continuous function $g:\RR\to\RR$,
\begin{equation}
\lim_{N\to\infty} \int_0^\infty g(s) P_N(s) ds = \int_0^\infty g(s) P(s) ds.
\end{equation}

The first question in convergence of probability measures is the problem of tightness.

\begin{lem}
The sequence of probability measures $\{ P_N(s) ds \}$ is tight on $\RR$. That is, for every $\epsilon>0$ there is a $K>0$ such that for all $N$
\begin{equation}
	\int_{|s|>K} P_N(s) ds < \epsilon.
\end{equation}
\end{lem}

\begin{proof}
We have 
\begin{equation}
\begin{split}
\int_{|s|>K} P_N(s) ds
& = \frac1N \#\{ j \leq N : s_j \geq K \} \\
& \leq \frac1N \sum_{j=1}^N \frac{s_j}{K} \chi_{[K,\infty)}(s_j) 
\leq \frac1N \sum_{j=1}^N \frac{s_j}{K} \\
& = \frac1K \sum_{j=1}^N (\xi_{j+1}-\xi_j) = \frac 1K.
\end{split}
\end{equation}
\end{proof}

Denote by $E_N(k,L)$ the probability of finding $k$ elements in the randomly shifted interval $[x_0,x_0+\tfrac{L}{N})$, i.e.,
\begin{equation}\label{EKL}
	E_N(k,L):=\meas\left\{ x_0\in\TT : 	S_N(\ell)=k \right\}.
\end{equation}
The following theorem explains the relation between $P(s)$ and the probability $E(0,L)$.

\begin{thm}\label{EP}
Given a probability density $P(s)$, 
the following statements are equivalent.
\begin{enumerate}
	\item[(i)] $P_N(s) \to_w P(s)$.
	\item[(ii)] $\lim_{N\to\infty} E_N(0,L) = E(0,L)$ for all $L>0$, where $E(0,L)$ is defined by
\end{enumerate}
\begin{equation}
\frac{d^2 E(0,L)}{dL^2} = P(L) , \qquad \lim_{L\to 0}E(0,L) = 1,  \qquad \lim_{L\to \infty} \frac{dE(0,L)}{dL} = 0.	
\end{equation}
\end{thm}

\begin{proof}
We have
\begin{equation}
\begin{split}
E_N(0,L) & =  \meas\big\{ x_0\in \TT^2 : \#\{j : \xi_j\in[x_0,x_0+\tfrac{L}{N})+\ZZ\}=0\big\} \\
& = \sum_{j=1}^N \meas\big\{ x_0\in [\xi_j,\xi_{j+1}) : \#\{j : \xi_j\in[x_0,x_0+\tfrac{L}{N})+\ZZ\}=0\big\} \\
& = \sum_{j=1}^N \left(\xi_{j+1}-\xi_j - \frac{L}{N}\right)
\chi_{[L,\infty)}(N(\xi_{j+1}-\xi_j)) \\
& = \sum_{j=1}^N \left(\xi_{j+1}-\xi_j\right)
- \sum_{j=1}^N \left(\xi_{j+1}-\xi_j\right) \chi_{[0,L)}(N(\xi_{j+1}-\xi_j)) \\
& - \frac{L}{N} \sum_{j=1}^N 
\chi_{[L,\infty)}(N(\xi_{j+1}-\xi_j))\\
& = 1 - \frac1N \sum_{j=1}^N g(s_j)
\end{split}
\end{equation}
where 
\begin{equation}
	g(x)= \max\{ 0 , x , L \} 
\end{equation}
is a bounded continuous function.

``(i)$\Rightarrow$(ii).'' With the above choice of test function $g$, (i) implies
\begin{equation}
	\lim_{N\to\infty} E_N(0,L) = F(L):=1 - \int_0^L s P(s) ds - L \int_L^\infty P(s) ds .
\end{equation}
Now 
\begin{equation}
\frac{dF(L)}{dL} = - \int_L^\infty P(s) ds , 	\qquad \frac{d^2 F(L)}{dL^2} = P(L) , 
\end{equation}
and 
\begin{equation}
\lim_{L\to 0}F(L) = 1,  \qquad \lim_{L\to\infty} \frac{dF(L)}{dL} = 0.
\end{equation}

``(ii)$\Rightarrow$(i).'' Since the sequence of probability measures $P_N(s)$ is tight, it is relatively compact by the Helly-Prokhorov Theorem (also often called Helly's Theorem). That is, every subsequence of $N$ contains a convergent subsequence $N_i$ for which $P_{N_i}(s)\to_w P(s)$ as $i\to\infty$. This implies (recall the first part of the proof) that $E_{N_i}(0,L) \to E(0,L)$ for all $L>0$. Hence every convergent subsequence has the limit $E(0,L)$, and thus every subsequence convergences.
\end{proof}

\subsection{Independent random variables}

In order to understand which statistical behaviour we should expect for the deterministic sequences we will study later, let us assume the vector $\xi=(\xi_1,\ldots,\xi_N)$ is a uniformly distributed random vector on $\TT^N$ with respect to Lebesgue measure $d\xi=d\xi_1\cdots d\xi_N$. (This means the $\xi_j$ are independent uniformly distributed random variables.) We can ignore the issue of ordering the $\xi_j$ here because of the symmetry of the measure $dx$ under permutation of coordinates. Expectation values and associated probabilities of a random variable $X=X(\xi)$ will be defined as
\begin{equation}
\EE X = \int_{\TT^N} X d\xi,
\end{equation}
\begin{equation}
\Prob(X>R) = \meas\{ \xi\in\TT^N : X>R \}.
\end{equation}

\begin{thm}
There is a constant $C>0$ such that, for all $\epsilon>0$, $N$, $L$, 
\begin{equation}
\Prob( |R_N^2(L,\psi)-L\widehat\psi(0)|>\epsilon ) \leq C \frac{L}{\epsilon^2 N}.
\end{equation}
\end{thm}

\begin{proof}
First of all, we have for the expectation (the $n=0$ term in \eqref{R22})
\begin{equation}
\EE R_N^2(L,\psi) = \frac{L(N-1)}{N} \widehat\psi(0)
= L \widehat\psi(0) (1+O(N^{-1})).
\end{equation}
Secondly, for the variance of $R_N^2(L,\psi)$,
\begin{multline}
\EE|R_N^2(L,\psi)-\EE R_N^2(L,\psi)|^2 \\
=\frac{L^2}{N^4} \sum_{\substack{i\neq j\\i'\neq j'}} \sum_{n,n'\in\ZZ} \widehat\psi\bigg(\frac{L n}{N}\bigg)\widehat\psi\bigg(\frac{L n'}{N}\bigg) \EE \big[e\big(n(\xi_i-\xi_j)-n'(\xi_{i'}-\xi_{j'})\big)\big].
\end{multline}
Now
\begin{equation}
\EE \big[e\big(n(\xi_i-\xi_j)-n'(\xi_{i'}-\xi_{j'})\big)\big]=
\begin{cases}
1 & \text{if $n=n'$, $i=i'$, $j=j'$} \\
&   \text{or if $n=-n'$, $i=j'$, $j=i'$} \\
0 & \text{otherwise.}
\end{cases}
\end{equation}
This implies that
\begin{equation}
\EE|R_N^2(L,\psi)-\EE R_N^2(L,\psi)|^2 
=\frac{L^2}{N^4} O(N^3/L) = O\big(\frac{L}{N}\big).
\end{equation}
\end{proof}

The above theorem implies that for a ``generic'' choice of the triangular array \eqref{triarray}, we have
\begin{equation}
R_N^2(L,\psi)=L\widehat\psi(0)+o(1)
\end{equation}
in the limit $N\to\infty$, $\ell=L/N\to 0$. This implies for the variance
\begin{equation}
\Sigma_N^2(L) = L +o(L) 
\end{equation}
almost surely in the above limit.

Using standard techniques from probability theory, one can extend these results on the variance to the full distribution of a generic realization of the random sequence in a small randomly shifted interval. There are two scaling regimes. 

{\bf Regime I (Central Limit Theorem):} In the limit $L\to\infty$, $N\to\infty$, $\ell=L/N\to 0$ we have
\begin{equation}
\meas\left\{ x_0\in\TT : 	\frac{S_N(\ell)-L}{\sqrt{\Sigma_N^2(L)}} > R \right\} \to 
\frac{1}{\sqrt{2\pi}} \int_R^\infty \e^{-t^2/2} dt 
\end{equation}
almost surely.

{\bf Regime II (Poisson Limit Theorem):} For $L$ fixed, $N\to\infty$, we have
\begin{equation}\label{poisson}
E_N(k,L) \to \frac{L^k}{k!} \e^{-L}. 
\end{equation}
almost surely.

\section{$m\alpha$ mod one}\label{secMA}

We will now consider the statistical properties of the sequence given by the fractional parts of $m\alpha$, $m=1,2,3,\ldots$ for some $\alpha$. This problem was studied by Berry-Tabor, Pandey et al., Bleher, Mazel-Sinai and Greenman using continued fractions (see \cite{npoint} for detailed references). In particular, it is a classical result that there are at most three distinct values for the gaps occurring in $m\alpha$ mod 1 which already indicates a rather non-generic behavior of the sequence, see e.g. \cite{Slater67}. 

Here we will use the approach introduced in \cite{npoint} that has the advantage of avoiding continued fractions and thus allowing higher-dimensional generalizations, such as the analysis of the distribution of linear forms modulo one. It is also very close to the work of Elkies and McMullen on $\sqrt m$ mod 1 which we will discuss in the next section.

We will be interested in the regime where $L=N\ell$ is fixed (Poisson scaling regime). The number \eqref{SN} of elements in an interval of size $\ell$ and centered at $x_0$ is then
\begin{equation}
\begin{split}
S_N(\ell) & = \sum_{m=1}^N \sum_{n\in\ZZ} \chi\bigg(\frac{N}{L} (m\alpha+n-x_0)\bigg) \\
& = \sum_{(m,n)\in\ZZ^2} \chi_{(0,1]}\bigg( \frac{m}{N} \bigg) \chi_{[-L/2,L/2]} \big(N (m\alpha+n-x_0)\big) \\
& = \sum_{(m,n)\in\ZZ^2} \psi\bigg( (m,n-x_0) \begin{pmatrix} 1 & \alpha \\ 0 & 1 \end{pmatrix} \begin{pmatrix} N^{-1} & 0 \\ 0 & N \end{pmatrix} \bigg)
\end{split}
\end{equation}
where $\chi_I$ denotes the characteristic function of the interval $I\subset\RR$ and 
\begin{equation}
	\psi(x,y) = \chi_{(0,1]}(x) \chi_{[-L/2,L/2]}(y)
\end{equation}
is the characteristic function of a rectangle.

Define the Lie Group $G$ by the semidirect product $\SL(2,\RR)\ltimes\RR^2$ with multiplication law
\begin{equation}
	(M,\xi)(M',\xi')=(MM',\xi M'+\xi'),
\end{equation}
where $\xi,\xi'\in\RR^2$ are viewed as row vectors.
This group has the matrix representation
\begin{equation}
	(M,\xi) \mapsto \begin{pmatrix} M & 0 \\ \xi & 1 \end{pmatrix} \in \SL(3,\RR).
\end{equation}
The function 
\begin{equation}\label{Fdef}
	F(M,\xi) = \sum_{m\in\ZZ^2} \psi(m M +\xi)
\end{equation}
defines a function on $G$. Note that, with $\psi$ as above, the sum in \eqref{Fdef} is always finite, and hence $F$ is a piecewise constant function. Furthermore,
\begin{equation}
	S_N(\ell) = F(M,\xi) 
\end{equation}
for the special choice
\begin{equation}
	M= \begin{pmatrix} 1 & \alpha \\ 0 & 1 \end{pmatrix} \begin{pmatrix} N^{-1} & 0 \\ 0 & N \end{pmatrix}, \qquad
	\xi=(0,-x_0) M.
\end{equation}

The crucial observation is now that $F$ is left-invariant under the discrete subgroup $\Gamma=\SL(2,\RR)\ltimes\RR^2$, and hence $F$ may be viewed as a piecewise constant function on the homogeneous space $\Gamma\backslash G$.

\begin{prop}
$F(\hat\gamma g)=F(g)$ for all $\hat\gamma\in\Gamma$.
\end{prop}

\begin{proof}
We have the decomposition
\begin{equation}
	\hat\gamma=(\gamma,n)=(\gamma,0)(1,n)
\end{equation}
for some $\gamma\in\SL(2,\ZZ)$, $n\in\ZZ^2$. It is therefore sufficient to check the statement for elements of the form $(\gamma,0)$ and $(1,n)$ separately.
We have
\begin{equation}
\begin{split}
	F((1,n)(M,\xi)) & = F(M,nM+\xi) \\
	& = \sum_{m\in\ZZ^2} \psi((m+n) M +\xi) \\
	& = \sum_{m\in\ZZ^2} \psi(m M +\xi)\\
	& = F(M,\xi)
\end{split}
\end{equation}
which proves one case, and
\begin{equation}
\begin{split}
	F((\gamma,0)(M,\xi)) & = F(\gamma M,\xi) \\
	& =  \sum_{m\in\ZZ^2} \psi(m \gamma M +\xi) \\
	& = \sum_{m\in\ZZ^2} \psi(m M +\xi) \\
	& = F(M,\xi)
\end{split}
\end{equation}
since $\gamma\ZZ^2=\ZZ^2$. 
\end{proof}

Alternatively, $F$ may be expressed as
\begin{equation}\label{Falt}
	F(g)=\sum_{\hat\gamma\in\pi(\Gamma)\backslash\Gamma} \psi(\pi(\hat\gamma g))
\end{equation}
with the projection 
\begin{equation}
	\begin{matrix} \pi: & G & \to & \RR^2 \\ & (M,\xi) & \mapsto & \xi . \end{matrix}
\end{equation}
From \eqref{Falt} the invariance under $\Gamma$ is directly evident.

\subsection{Geometry of $\Gamma\backslash G$}

\begin{figure}
\includegraphics[width=0.5\textwidth]{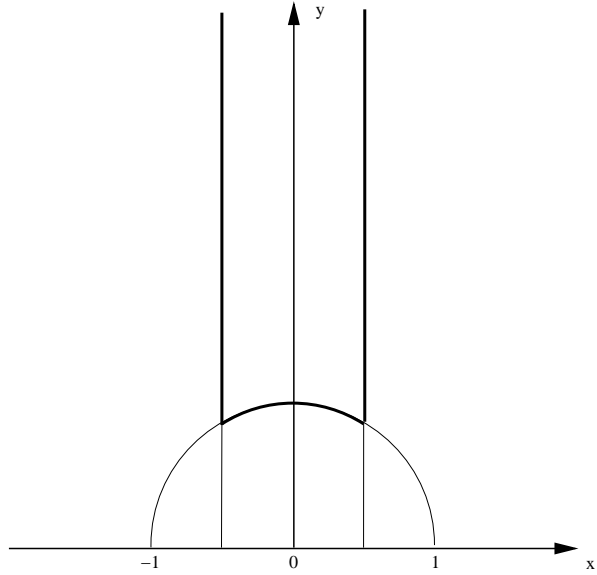}
\caption{Fundamental domain of the modular group $\SL(2,\ZZ)$ in the complex upper half plane.\label{mod}}
\end{figure}

The aim is to find a good coordinate system for $G$. Since parametrizing $\RR^2$ is obvious, we need to mainly worry about $\SL(2,\RR)$. The Iwasawa decomposition of an element $M\in\SL(2,\RR)$ is
\begin{equation}
	M=\begin{pmatrix} 1 & u \\ 0 & 1 \end{pmatrix}\begin{pmatrix} v^{1/2} & 0 \\ 0 & v^{-1/2} \end{pmatrix}
	\begin{pmatrix} \cos(\phi/2) & \sin(\phi/2) \\ -\sin(\phi/2)  & \cos(\phi/2) \end{pmatrix}
\end{equation}
where $\tau=u+\i v \in\HH:=\{ \tau \in\CC: \Im\tau>0\}$ (the complex upper halfplane) and $\phi\in[0,4\pi)$.
This yields a 1-1 map $\SL(2,\RR)\to\HH\times [0,4\pi)$. Left-multiplication becomes now an action of $\SL(2,\RR)$ on $\HH\times [0,4\pi)$ given by the formula
\begin{equation}
	\begin{pmatrix} a & b \\ c & d \end{pmatrix} = \left( \frac{a\tau+b}{c\tau+d}, \phi-2\arg(c\tau+d) \right)
\end{equation}
(this can be checked by a straightforward calculation).
The fractional linear transformation of the $\tau$ component defines an (orientation preserving) isometry with respect to the Riemannian line element
\begin{equation}
	ds^2=\frac{du^2+dv^2}{v^2}
\end{equation}
and the transformation property of $\phi$ is identical to the direction of a tangent vector at $\tau\in\HH$.
Thus the group $\PSL(2,\RR):=\SL(2,\RR)/\{\pm1\}\simeq \HH\times [0,2\pi)$ can be identified with the unit tangent bundle $\T^1\HH$ of $\HH$. Similarly, $\SL(2,\ZZ)\backslash\SL(2,\RR)\simeq \PSL(2,\ZZ)\backslash\PSL(2,\RR)$ can be identified with the unit tangent bundle of the modular surface $\SL(2,\ZZ)\backslash\HH$. A fundamental domain $\scrF$ for the action of $\SL(2,\ZZ)$ on $\HH$ is shown in Figure \ref{mod}. We have
\begin{equation}
\begin{split}
	\scrF & =\{ \tau\in\HH : |\tau| > 1,\; |\Re \tau|< 1/2 \} \\
	& \cup \{ \tau\in\HH : |\tau| \geq 1,\; \Re \tau=-1/2\} \\
	& \cup \{ \tau\in\HH : |\tau| = 1,\; -1/2\leq \Re \tau \leq 0 \}.
\end{split}
\end{equation}
Note that the modular surface is not compact, there is one cusp at $\i\infty$. It has however finite measure with respect to the Riemannian volume $v^{-2} du\,dv$.

In order to understand the geometry of all of $\GamG$, write
\begin{equation}
	g=(1,\xi)(M,0)
\end{equation}
which gives a particular parametrization in terms of $\RR^2$ and $\SL(2,\RR)$. Since $\Gamma$ contains the subgroup $1\ltimes \ZZ^2$, $\xi$ can be parametrized by $\TT^2=\ZZ^2\backslash\RR^2\simeq [0,1)^2$.
This concludes our analysis: we have found a 1-1 parametrization of $G$ in terms of
\begin{equation}
	\T^1(\SL(2,\RR)\backslash\HH)\times \TT^2 \simeq \scrF \times [0,2\pi) \times [0,1)^2 .
\end{equation}
That is, $\GamG$ is a (non-trivial) bundle over $\T^1(\SL(2,\ZZ)\backslash\HH)$ with fibre $\TT^2$.

\subsection{Dynamics on $\GamG$}

Consider the one-parameter subgroup $\Phi^\RR:=\{ \Phi^t \}_{t\in\RR}$ where
\begin{equation}
	\Phi^t = \left( \begin{pmatrix} \e^{-t/2} & 0 \\ 0 & \e^{t/2} \end{pmatrix}, 0 \right).
\end{equation}
$\Phi^\RR$ defines a flow on $\GamG$ by right multiplication,
\begin{equation}
	\Gamma g \mapsto \Gamma g \Phi^t.
\end{equation}
The remarkable observation is that our object of interest, $S_N(\ell)$, is related to a function $F$ on $\GamG$ evaluated along an orbit of this flow:
\begin{equation}
	S_N(\ell) = F(g_0 \Phi^t)
\end{equation}
with $t=2\log N$ and initial condition 
\begin{equation}
	g_0 = \left( \begin{pmatrix} 1 & \alpha \\ 0 & 1 \end{pmatrix}, (0,-x_0) \right) .
\end{equation}
Let us define
\begin{equation}
	n_-(\alpha,y) = \left( \begin{pmatrix} 1 & \alpha \\ 0 & 1 \end{pmatrix}, (0,y) \right) .
\end{equation}
The subgroup $H=\{n_-(\alpha,y)\}_{(\alpha,y)\in\RR^2}$ is abelian and isomorphic to $\RR^2$. Notice that 
\begin{equation}
	\Gamma\cap H = \{n_-(\alpha,y)\}_{(\alpha,y)\in\ZZ^2}
\end{equation}
is a subgroup of $H$ isomorphic to $\ZZ^2$. Therefore, for every fixed $t$, the set
\begin{equation}
\Gamma\backslash\Gamma H \Phi^t
\end{equation}
describes a torus $\simeq\TT^2$ embedded in $\GamG$; $t$ parametrizes a continuous family of such tori.

We will now show that $H$ parametrizes the unstable directions of the flow $\Phi^t$. We employ the following parametrization of $G$. Write
\begin{equation}
	g= n_-(\alpha,y) \Phi^s n_+(\beta,x), 
\end{equation}
where
\begin{equation}
	n_+(\beta,x) = \left( \begin{pmatrix} 1 & 0 \\ \beta & 1 \end{pmatrix}, (x,0) \right) .
\end{equation}
We will write for short $g=(\alpha,y,s,\beta,x)$. The advantage of these coordinates is that the time evolution under $\Phi^t$ can be worked out very simply. We have the relation
\begin{equation}\label{phiphi}
	(\alpha,y,s,\beta,x)\Phi^t=\Phi^t (\e^t \alpha,\e^{t/2} y,s,\e^{-t}\beta,\e^{-t/2}x).
\end{equation}
Distances on $\GamG$ are measured by a left-$G$-invariant (since $\Gamma$ acts on the left) Riemannian metric $d(g,g')$ on $G$. If $g=(\alpha,y,0,0,0)$ and $g'=(\alpha',y',0,0,0)$ are two initially close points, we have under the flow $\Phi^t$ (use the above formula and left-invariance of the metric)
\begin{equation}
\begin{split}
	d(g\Phi^t,g'\Phi^t) & = d((\e^t(\alpha-\alpha'),\e^{t/2}(y-y'),0,0,0),(0,0,0,0,0)) \\
	& \approx (\e^{2t} |\alpha-\alpha'|^2+\e^{t}|y-y'|^2)^{1/2}.
\end{split}
\end{equation}
Hence $(\alpha,y)$ describe exponentially unstable directions of the flow, and by the same argument it is easy to see that $(\beta,y)$ are the exponentially stable directions and $s$ is of course the neutral flow direction. In particular we have the bound
\begin{equation}\label{dell}
	d((\alpha,y,s,\beta,x)\Phi^t,(\alpha,y,0,0,0)\Phi^t) =O(|s|+|\beta|\e^{-t} + |x|\e^{-t/2})
\end{equation}
for $s,\beta,x$ bounded and $t>0$. This follows directly from \eqref{phiphi}.

\subsection{Mixing and uniform distribution}

Recall that we are interested in the behaviour of the distribution of $S_N(\ell)$ for $x_0$ random and $N$ large. At this point it will be convenient to also take $\alpha$ to be random, say, uniformly distributed in the interval $[a,b]$. We will see later that for fixed $\alpha$ there is no universal limiting distribution (an observation that is well known and related to the three gap theorem \cite{Slater67}).

We will use equidistribution on $\GamG$ to prove the following limit theorem, which asserts a limiting distribution different from Poissonian, cf. \eqref{poisson}.

We will use the notation $\overline g=\Gamma g$.

\begin{thm}\label{limthm}
For any $L>0$,
\begin{equation}
	\lim_{N\to\infty} \frac{1}{b-a} \meas\{ (\alpha,x_0)\in[a,b]\times[0,1] : S_N(\ell)=k\} = E(k,L),
\end{equation}
where
\begin{equation}\label{EF}
	E(k,L)=\frac{1}{\mu(\GamG)} \mu( \overline g\in\GamG : F(g)=k ).
\end{equation}
\end{thm}

Here $F$ is the function defined in \eqref{Fdef}, and $\mu$ the Haar measure on $G$. An explicit formula for $d\mu$ in the Iwasawa coordinates is
\begin{equation}
	d\mu = \frac{du\,dv}{v^2} \,d\phi\,dx\,dy .
\end{equation}
It is possible to derive more explicit formulas for $E(k,L)$ from \eqref{EF}, but this requires some involved calculations which we will not pursue her. See \cite{Strombersson05}, Section 8, for details.

The key to the proof is the following equidistribution theorem.

\begin{thm}\label{udi}
For any bounded, piecewise continuous\footnote{i.e. the discontinuities are contained in a set of $\mu$ measure zero.} $f:\GamG\to\RR$
\begin{equation}
	\lim_{t\to\infty} \frac{1}{b-a} \int_a^b\int_0^1 f(n_-(\alpha,y)\Phi^t) d\alpha\,dy
	=\frac{1}{\mu(\GamG)} \int_{\GamG} f d\mu.
\end{equation}
\end{thm}

\begin{proof}
It is well known that the flow $\Phi^t$ is mixing,\footnote{This is guaranteed by a general theorem by Moore for semisimple Lie groups, which can be extended to the non-semisimple $G$ considered here, cf. \cite{Kleinbock99}.} that is
for any $f,h\in\L^2(\GamG)$
\begin{equation}
		\lim_{t\to\infty} \int_{\GamG} f(g\Phi^t) h(g) d\mu
	=\frac{1}{\mu(\GamG)} \int_{\GamG} f d\mu \int_{\GamG} h d\mu.
\end{equation}
Take $f$ to be continuous an of compact support, and $h$ the characteristic function of the set
\begin{equation}
	S_\epsilon=\Gamma\{ (\alpha,y,s,\beta,y): \alpha\in[a,b],\,y\in[0,1], s,\beta,x\in[-\epsilon,\epsilon]\},
\end{equation}
which forms an $\epsilon$-neighbourhood of the embedded closed torus $S_0$. 
By the uniform continuity of $f$ and \eqref{dell}, given any $\delta>0$ there is an $\epsilon>0$ such that
\begin{equation}
	\sup_{\substack{g\in S_\epsilon \\ t>0}} |f(g\Phi^t)-f(n_-(\alpha,y)\Phi^t)| < \delta.
\end{equation}
Haar measure in the local coordinates $(\alpha,y,s,\beta,y)$ reads (up to normalization) 
\begin{equation}
	d\mu=\e^{3s/2} ds\,d\alpha\,d\beta\,dx\,dy.
\end{equation}
We conclude that
\begin{equation}
	\liminf_{t\to\infty} \frac{1}{b-a} \int_a^b\int_0^1 f(n_-(\alpha,y)\Phi^t) d\alpha\,dy
	=\frac{1}{\mu(\GamG)} \int_{\GamG} f d\mu +O(\delta)
\end{equation}
and
\begin{equation}
	\limsup_{t\to\infty} \frac{1}{b-a} \int_a^b\int_0^1 f(n_-(\alpha,y)\Phi^t) d\alpha\,dy
	=\frac{1}{\mu(\GamG)} \int_{\GamG} f d\mu + O(\delta),
\end{equation}
where the implied constants are independent of $\epsilon$. This works for any $\delta>0$, and hence the limit must exist and equal $\frac{1}{\mu(\GamG)} \int_{\GamG} f d\mu$. 

To extend the statement of the theorem to bounded continuous functions, we observe that it holds (trivially) for constant $f$, and therefore also for continuous functions $f$ that are constant outside some compact set. 

Let $f$ be a bounded piecewise continuous function. Given any $\epsilon>0$ we can find continuous functions $f_\pm$, constant outside some constant set, such that
\begin{equation}
	f_- \leq f \leq f_+
\end{equation}
and
\begin{equation}
\frac{1}{\mu(\GamG)}	\int_{\GamG} (f_+-f_-) d\mu < \epsilon.
\end{equation}
This implies
\begin{equation}
\begin{split}
	\liminf_{t\to\infty} \frac{1}{b-a} \int_a^b\int_0^1 f(n_-(\alpha,y)\Phi^t) d\alpha\,dy
	& \geq \liminf_{t\to\infty} \frac{1}{b-a} \int_a^b\int_0^1 
	f_-(n_-(\alpha,y)\Phi^t) d\alpha\,dy \\
	& = \frac{1}{\mu(\GamG)} \int_{\GamG} f_- d\mu \\
	& > \frac{1}{\mu(\GamG)} \int_{\GamG} f d\mu - 2\epsilon.
\end{split}
\end{equation}
The analogous argument shows
\begin{equation}
	\limsup_{t\to\infty} \frac{1}{b-a} \int_a^b\int_0^1 f(n_-(\alpha,y)\Phi^t) d\alpha\,dy
	< \frac{1}{\mu(\GamG)} \int_{\GamG} f d\mu + 2\epsilon.
\end{equation}
Taking $\epsilon>0$ arbitrarily small proves the theorem.
\end{proof}

\begin{remark}
An alternative proof of Theorem \ref{udi} follows from Ratner's theorem, since the subgroup $\{ n_-(\alpha,y)\}_{\alpha,y\in\RR}$ is generated by unipotent elements. We will get back to this later.
\end{remark}

\begin{proof}[Proof of Theorem \ref{limthm}]
Apply Theorem \ref{udi} to the characteristic function of the set of $\overline g\in\GamG$ for which $F(g)=k$ (to make sure the characteristic function is piecewise continuous, check that the set has a boundary of $\mu$ measure zero).
\end{proof}

\begin{remark}
As we had mentioned earlier, there is no limiting distribution as in Theorem \ref{limthm} if $\alpha$ is fixed, since there is no analog of the equidistribution result, Theorem \ref{udi}. One can show, however, that if $\alpha$ is irrational we have for any continuous, compactly supported function
\begin{equation}
	\int_{\TT} f(n_-(\alpha,y)\Phi^t) dy
	= \overline f\left(\begin{pmatrix} 1 & \alpha \\ 0 & 1 \end{pmatrix}
	\begin{pmatrix} \e^{-t/2} & 0 \\ 0 & \e^{t/2} \end{pmatrix}\right) +o(1) \qquad(t\to\infty)
\end{equation}
where $\overline f$ is a (non-constant!) continuous, compactly supported function on $\SL(2,\ZZ)\backslash\SL(2,\RR)$ defined by
\begin{equation}
	\overline f(M)= \int_{\TT^2} f( (1,\xi) (M,0)) d\xi .
\end{equation}
\end{remark}

\begin{remark}\label{rem2}
If one however fixes $y=-x_0\notin\QQ$ and keeps $\alpha$ random, Ratner's Theorem implies the following equidistribution result.
For any bounded piecewise continuous $f:\GamG\to\RR$
\begin{equation}\label{udi00}
	\lim_{t\to\infty} \int_{\TT} f(n_-(\alpha,y)\Phi^t) d\alpha
	=\frac{1}{\mu(\GamG)} \int_{\GamG} f d\mu.
\end{equation}
Hence the limiting distribution is universal (i.e. independent of $y$ as long as $y$ is irrational) and the same as for random $y$. Thus the probability of finding $k$ points in the interval $x_0-\ell/2,x_0+\ell/2)$ with {\em fixed} center $x_0\notin\QQ$ has the limiting distribution 
\begin{equation}
	\lim_{N\to\infty} \meas\{ \alpha\in\TT^2 : S_N(\ell)=k\} = E(k,L),
\end{equation}
the same as for {\em random} center. We will prove \eqref{udi} in Section \ref{secRatner}.
\end{remark}

\section{$\sqrt{m\alpha}$ mod one}

The problem of the statistics of $\sqrt{m\alpha}$ mod 1 has been understood by Elkies and McMullen \cite{Elkies04} in the case $\alpha=1$ (and in principle also for all other rational $\alpha$). The uniform distribution of $\sqrt{m\alpha}$ mod 1 may be shown by using the fact that $\sqrt{n+m}-\sqrt{n}\to 0$ for $n\to\infty$, $m$ fixed (we leave this as an exercise). As in the last section, the key idea is the reduce the problem to equidistribution on a homogeneous space. Lucky for us, this homogeneous space will turn out to be $\GamG$ with the same $G$, $\Gamma$ as encountered earlier.

We are as in the previous section interested in the ``Poisson scaling limit'', i.e. $L$ is fixed. Now (we swap $m$ and $n$ in our notation)
\begin{equation}
	S_N(\ell)= \sum_{n=1}^N \sum_{m\in\ZZ} \chi\left( \frac{N}{L} ( \sqrt{n\alpha} - x_0 +m) \right).
\end{equation}
The condition imposed on the summation can be re-written as
\begin{equation}
	\left(x_0-m-\frac{L}{2N}\right)^2 \leq n \alpha < \left(x_0-m+\frac{L}{2N}\right)^2
\end{equation}
which amounts to
\begin{equation}
	-\frac{L}{N} (x_0-m) \leq n\alpha -(x_0-m)^2 - \left(\frac{L}{2N}\right)^2 < \frac{L}{N} (x_0-m) .
\end{equation}
Notice also that
\begin{equation}
	|\sqrt{n\alpha}-(x_0-m)|\leq \frac{L}{2N}.
\end{equation}
This yields
\begin{multline}
	S_N(\ell)= \sum_{(m,n)\in\ZZ^2} \chi_{(0,1]}\left(\frac{x_0-m+O(L/2N)}{\sqrt{N\alpha}}\right)\\
	\chi_{[-L,L)}\left( \frac{N^{1/2}[n\alpha - (x_0 -m)^2-(L/2N)^2]}{N^{-1/2}(x_0-m)} \right) .
\end{multline}
A more convenient object would be
\begin{equation}
	\widetilde S_{N,\epsilon,\delta}(\ell)= \sum_{(m,n)\in\ZZ^2} \chi_{(-\epsilon,1+\epsilon]}\left(\frac{x_0-m}{\sqrt{N\alpha}}\right)
	\chi_{[-L,L)}\left( \frac{N^{1/2}[n\alpha - (x_0 -m)^2]+\delta}{N^{-1/2}(x_0-m)} \right) .
\end{equation}
For the right choices of $\epsilon$ (positive/negative) we obtain upper/lower bounds for $S_N(\ell)$ which would eventually allow us to infer the limiting distribution of $S_N(\ell)$ from $\widetilde S_{N,\epsilon,\delta}(\ell)$ by taking $\delta\to 0$, $\epsilon\to \pm 0$. We will ignore this technical point here and simply take
\begin{equation}\label{approx}
	S_N(\ell)\approx\widetilde S_{N,0,0}(\ell) = \sum_{(m,n)\in\ZZ^2} \chi_{(0,1]}\left(\frac{x_0-m}{\sqrt{N\alpha}}\right)
	\chi_{[-L,L)}\left( \frac{N^{1/2}[n\alpha - (x_0 -m)^2]}{N^{-1/2}(x_0-m)} \right) .
\end{equation}
The manipulations we will now perform on the r.h.s. of \eqref{approx} can be adapted step by step for more general values of $\delta,\epsilon\neq 0$ (recommended exercise). We will use the shorthand $\widetilde S_N(\ell):=\widetilde S_{N,0,0}(\ell)$ in the following.

\subsection{The case $\alpha=1$}

We have, after substituting $(m,n)\to(-m,-n)$,
\begin{equation}\label{approx1}
	\widetilde S_N(\ell) = \sum_{(m,n)\in\ZZ^2} \chi_{(0,1]}\left(\frac{x_0+m}{\sqrt{N}}\right)
	\chi_{(-L,L]}\left( \frac{N^{1/2}[n+(x_0 +m)^2]}{N^{-1/2}(x_0+m)} \right) ,
\end{equation}
an thus, after substituting $n\mapsto n+m^2$ in the sum over $n$,
\begin{equation}\label{approx2}
	\widetilde S_N(\ell) = \sum_{(m,n)\in\ZZ^2} \chi_{(0,1]}\left(\frac{x_0+m}{\sqrt{N}}\right)
	\chi_{(-L,L]}\left( \frac{N^{1/2}(n + x_0^2 + 2mx_0)}{N^{-1/2}(x_0+m)} \right) ,
\end{equation}
We will now show that, in analogy with the previous section, we can find a function $F:\GamG\to\RR$ of the form
\begin{equation}\label{FF}
	F(M,\xi) = \sum_{m\in\ZZ^2} \psi(mM+\xi)
\end{equation}
so that
\begin{equation}\label{Fg}
	\widetilde S_N(\ell) = F(g)
\end{equation}
for a suitable choice of $g\in G$ and a piecewise continuous $\psi:\RR^2\to\RR$ with compact support. To this end define
\begin{equation}
	\psi(x,y)=\chi_{(0,1]}(x) \chi_{(-L,L]}\left(\frac{y}{x}\right)
\end{equation}
(which indeed has compact support: it is the characteristic function of a triangle). Now
consider the one parameter subgroup $\{ n_1(x) \}_{x\in\RR}$ with
\begin{equation}
	n_1(x)= \left(\begin{pmatrix} 1 & 2x \\ 0 & 1 \end{pmatrix}, (x,x^2) \right) 
\end{equation}
(check that this indeed yields a one parameter group). 
Then the choice (set $t=\log N$, $x=x_0$)
\begin{equation}
\begin{split}
	(M,\xi) & = n_1(x_0)\Phi^t \\
	  & = \left(\begin{pmatrix} N^{-1/2} & 2x_0 N^{1/2} \\ 0 & N^{1/2} \end{pmatrix}, (N^{-1/2} x_0, N^{1/2} x_0^2) \right)
\end{split}
\end{equation}
yields 
\begin{equation}
	(m,n)M+\xi=\big( N^{-1/2} (x_0+m), N^{1/2} (2mx_0+n+x_0^2)\big) .
\end{equation}
Using this result in the definition \eqref{FF} then confirms the desired \eqref{Fg}.

We now follow the same steps as in the previous Section \ref{secMA} to derive the limiting distribution for $S_N(\ell)$ from equidistribution on $\GamG$. We first state the limit theorem.

\begin{thm}\label{limthm2}
For any $L>0$,
\begin{equation}
	\lim_{N\to\infty} \meas\{ x_0\in\TT : S_N(\ell)=k\} = E(k,L),
\end{equation}
where
\begin{equation}\label{EF2}
	E(k,L)=\frac{1}{\mu(\GamG)} \mu( \overline g\in\GamG : F(g)=k ),
\end{equation}
with $F$ as defined in \eqref{FF}.
\end{thm}

An explicit formula for $E(0,L)$ and the corresponding gap distribution $P(s)$ (recall Theorem \ref{EP}) is worked out in \cite{Elkies04}.

The relevant equidistribution theorem needed to prove Theorem \ref{limthm2} is the following. Note that $\Gamma \cap \{n_1(x)\}_{x\in\RR}= \{n_1(x)\}_{x\in\ZZ}$ and hence 
\begin{equation}
\Gamma\{n_1(x)\}_{x\in\TT}\Phi^t 
\end{equation}
represents a family (parametrized by $t$) of closed orbits embedded in $\GamG$.

\begin{thm}\label{udi2}
For any bounded piecewise continuous $f:\GamG\to\RR$
\begin{equation}
	\lim_{t\to\infty} \int_{\TT} f(n_1(x)\Phi^t) dx
	=\frac{1}{\mu(\GamG)} \int_{\GamG} f d\mu.
\end{equation}
\end{thm}

Since $n_1(x)$ generates a unipotent flow, Ratner's theorem can be employed. We will explain the general strategy of proof for statements of this type in Section \ref{secRatner}.

\subsection{Some heuristics in the case $\alpha=\sqrt 2$}

\begin{figure}
\begin{minipage}{0.8\textwidth}
\footnotesize
%% Created by Maple 8.00 (IBM INTEL NT)
%% Source Worksheet: sqrtnalpha.mws
%% Generated: Tue Apr 04 13:11:48 2006
\begin{maplegroup}
\end{maplegroup}
\begin{maplegroup}
\begin{mapleinput}
\mapleinline{active}{1d}{alpha:=sqrt(2); N:=6001;}{%
}
\end{mapleinput}

\mapleresult
\begin{maplelatex}
\mapleinline{inert}{2d}{alpha := 2^(1/2);}{%
\[
\alpha  := \sqrt{2}
\]
}
\end{maplelatex}

\begin{maplelatex}
\mapleinline{inert}{2d}{N := 6001;}{%
\[
N := 6001
\]
}
\end{maplelatex}

\end{maplegroup}
\begin{maplegroup}
\begin{mapleinput}
\mapleinline{active}{1d}{L:=sort([seq(evalf[12](frac(sqrt(n*alpha))), n=1..N)]):}{%
}
\end{mapleinput}

\end{maplegroup}
\begin{maplegroup}
\end{maplegroup}
\begin{maplegroup}
\begin{mapleinput}
\mapleinline{active}{1d}{alist:=seq(evalf[12](N*(L[i+1]-L[i])),i=1..N-1):}{%
}
\end{mapleinput}

\end{maplegroup}
\begin{maplegroup}
\begin{mapleinput}
\mapleinline{active}{1d}{data:=stats[transform,tallyinto['outliers']]([alist],[seq((i-1)*0.2..
i*0.2,i=0..35)]):}{%
}
\end{mapleinput}

\end{maplegroup}
\begin{maplegroup}
\begin{mapleinput}
\mapleinline{active}{1d}{outliers;}{%
}
\end{mapleinput}

\mapleresult
\begin{maplelatex}
\mapleinline{inert}{2d}{[7.0547245915, 7.0674227075, 7.1105849000, 7.1693268887,
7.2093775627, 7.3219323187, 7.3381866273, 7.4195061783, 7.5000233956,
7.6451419780, 7.7497418084, 7.9388213164, 8.0221013941, 8.1512135092,
8.4582030656];}{%
\maplemultiline{
[7.0547245915, \,7.0674227075, \,7.1105849000, \,7.1693268887, \,
7.2093775627,  \\
7.3219323187, \,7.3381866273, \,7.4195061783, \,7.5000233956, \,
7.6451419780,  \\
7.7497418084, \,7.9388213164, \,8.0221013941, \,8.1512135092, \,
8.4582030656] }
}
\end{maplelatex}

\end{maplegroup}
\begin{maplegroup}
\begin{mapleinput}
\mapleinline{active}{1d}{data1:=stats[transform,scaleweight[1/nops([alist])]](data):}{%
}
\end{mapleinput}

\end{maplegroup}
\begin{maplegroup}
\end{maplegroup}
\begin{maplegroup}
\begin{mapleinput}
\mapleinline{active}{1d}{g1:=stats[statplots,histogram](data1):}{%
}
\end{mapleinput}

\end{maplegroup}
\begin{maplegroup}
\begin{mapleinput}
\mapleinline{active}{1d}{g2:=plot(exp(-s), s=0..6):}{%
}
\end{mapleinput}

\end{maplegroup}
\begin{maplegroup}
\begin{mapleinput}
\mapleinline{active}{1d}{plots[display](g1,g2);}{%
}
\end{mapleinput}

\mapleresult
\begin{center}
\mapleplot{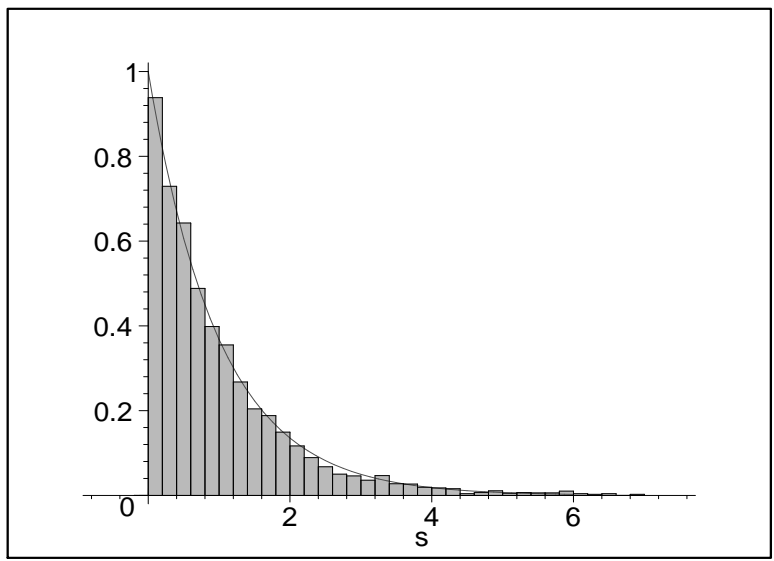}
\end{center}

\end{maplegroup}
\begin{maplegroup}
\end{maplegroup}
\begin{maplegroup}
\end{maplegroup}
\begin{maplegroup}
\end{maplegroup}
%% End of Maple 8.00 Output
\end{minipage}
\caption{Maple worksheet for calculating the gap distribution of the fractional parts of $\sqrt{m\sqrt2}$, $m=1,\ldots,6001$. \label{maple}}
\end{figure}

We return to generic $\alpha$, such as $\alpha=\sqrt 2$, and rewrite $\widetilde S_N(\ell)$ as
\begin{equation}
	\widetilde S_N(\ell)=\sum_{(m,n)\in\ZZ^2} \chi_{(0,1]}\left(\frac{x_0+m}{M}\right)
	\chi_{(-L,L]}\left( \frac{M[\alpha^{-1} (x_0 +m)^2+n]}{M^{-1}(x_0+m)} \right) 
\end{equation}
where $M=\sqrt{N\alpha}$. For $x_0\in[0,1]$ we can ignore terms of the form $x_0/M$, 
\begin{equation}
	\widetilde S_N(\ell)\approx \sum_{(m,n)\in\ZZ^2} \chi_{(0,1]}\left(\frac{m}{M}\right)
	\chi_{(-L,L]}\left( \frac{M[\alpha^{-1} (x_0 +m)^2+n]}{M^{-1}m} \right) .
\end{equation}
Now note that for most values of $m$, we have $m/M\asymp 1$, and it is natural to assume that, for random $x_0$,  the probability of finding $k$ elements of the set 
\begin{equation}
	\{\alpha^{-1} (x_0 +m)^2:m=1,\ldots,M\} + \ZZ
\end{equation}
in an interval of size $1/M$ around the origin is given by the Poisson distribution (we must assume here that $\alpha$ is badly approximable by rationals, e.g. $\alpha=\sqrt 2$ would be a good choice). Hence we may assert that the limiting distribution of $S_N(\ell)$ is the same as that of the random variable
\begin{equation}
	X=\sum_{(m,n)\in\ZZ^2} \chi_{(0,1]}\left(\frac{m}{M}\right)
	\chi_{(-L,L]}\left( \frac{M(\eta_m+n)}{M^{-1}m} \right) 
\end{equation}
where $\eta_m$ are independent uniformly distributed random variables on $[-1/2,1/2)$. With this choice of interval the only contribution comes from the $n=0$ term (assume $M\gg L$), so
\begin{equation}
	X=\sum_{m=1}^M X_m 
\end{equation}
where
\begin{equation}
	X_m=\chi_{(-L,L]}\left( \frac{M^2 \eta_m}{m} \right)
\end{equation}
is a sequence of independent random variables with $k$th moment
\begin{equation}
	\EE X_m^k = \int_{-1/2}^{1/2} \chi_{(-L,L]}\left( \frac{M^2 \eta_m}{m} \right) d\eta_m =\frac{2L m}{M^2} ,
\end{equation}
and hence 
\begin{equation}
	\EE(\e^{\i t X_m}-1)= \frac{2L m}{M^2} (\e^{\i t}-1).
\end{equation}
The characteristic function of the random variable $X$ is therefore
\begin{equation}
\begin{split}
	\EE \e^{\i t X} & = \prod_{m=1}^M \bigg[ 1 + \frac{2L m}{M^2} (\e^{\i t}-1) \bigg] \\ 
	& = \exp \bigg\{ \sum_{m=1}^M \log\bigg[ 1 + \frac{2L m}{M^2} (\e^{\i t}-1) \bigg] \bigg\}\\
	& = \exp \bigg\{ \sum_{m=1}^M \bigg[ \frac{2L m}{M^2} (\e^{\i t}-1) +O\bigg(\frac{m^2}{M^4}\bigg)\bigg]\bigg\} \\
	& = \exp\bigg[ L (\e^{\i t}-1) +O\bigg(\frac{1}{M}\bigg) \bigg] .
\end{split}
\end{equation}
The expression $\e^{L (\e^{\i t}-1)}$ is the characteristic function of the Poisson law
\begin{equation}
E(k,L)=	\frac{L^k}{k!} \e^{-L} .
\end{equation}
Hence this should be our prediction for the limiting distribution of $S_N(\ell)$, which in turn implies that we expect the exponential distribution for gaps in $\sqrt{m\alpha}$ mod 1. This is in good agreement with our Maple experiment, Figure \ref{maple}.

\section{Ratner's theorem}\label{secRatner}

An excellent introduction to Ratner's theory is Dave W. Morris' recent textbook \cite{Morris05}.
Let $G$ be a Lie group (e.g. $\SL(2,\RR)\times\RR^2$) and $\Gamma$ be a discrete subgroup (e.g. $\SL(2,\ZZ)\times\ZZ^2$). It is at this point not necessary to assume that $\Gamma$ is a lattice in $G$, i.e., that $\GamG$ has finite volume with respect to Haar measure $\mu$ on $G$. Ratner's measure classification theorem gives a complete geometric description of all measures that are invariant and ergodic under the a unipotent one parameter subgroup $U$ (or, more generally, invariant and ergodic under a subgroup generated by unipotent subgroups) acting on $\GamG$ by right multiplication. Examples of unipotent subgroups that appeared in the previous sections are $\{n_-(\alpha,0)\}_{\alpha\in\RR}$, $\{n_-(0,y)\}_{\alpha\in\RR}$ and $\{n_1(x)\}_{x\in\RR}$.

\begin{thm}[Ratner's theorem]
	Let $\nu$ be an ergodic, $U$-invariant probability measure on $\GamG$. Then there is a closed, connected subgroup $H\subset G$, and a point $\overline g\in\GamG$ such that
\begin{enumerate}
	\item $\nu$ is $H$-invariant,
	\item $\nu$ is supported on the orbit $\overline g H$.
\end{enumerate}
\end{thm}

\begin{remark}
Let $g\in G$ be a representative of the coset $\overline g=\Gamma g$, and define the subgroup $\Gamma_H=(g^{-1}\Gamma g)\cap H$. Then the orbit $\overline g H$ may be identified with the homogeneous space $\Gamma_H\backslash H$ and $\nu$ with the Haar measure on $H$. Furthermore one can deduce (since $\nu$ is a probability measure) that $\Gamma_H$ is a lattice in $H$, i.e., $\nu(\overline g H)<\infty$, and that the orbit $\overline g H$ is closed in $\GamG$. 
\end{remark}

In simple words, measures $\nu$ invariant and ergodic under unipotent subgroups are supported on nice embedded closed subvarieties, of which there can be only countably many (modulo translations of course). We will now discuss two corollaries of Ratner's theorem that are relevant to the equidistribution theorems discussed earlier.

\subsection{Limit distributions of translates}

The following is special case of Shah's extremely useful theorem, Theorem 1.4 in \cite{Shah96}.

\begin{thm}
Suppose $G$ contains a Lie subgroup $H$ isomorphic to $\SL(2,\RR)$ (we denote the corresponding embedding by $\varphi:\SL(2,\RR)\to G$), such that the set $\Gamma\backslash\Gamma H$ is dense in $\GamG$. 
Then, for any bounded, piecewise continuous $f:\GamG\to \RR$ and any piecewise continuous $h:\RR\to\RR$ with compact support
\begin{equation}
	\lim_{t\to\infty} \int_{\RR} f\left(\varphi\left(
	\begin{pmatrix} 1 & x \\ 0 & 1\end{pmatrix}
	\begin{pmatrix} \e^{-t/2} & 0 \\ 0 & \e^{t/2} \end{pmatrix}
	\right)\right) h(x) dx
	= \frac{1}{\mu(\GamG)} \int_{\GamG} f d\mu \, \int_{\RR} h(x) dx
\end{equation}
where $\mu$ is the Haar measure of $G$.
\end{thm}

The general strategy of proof for statements of the above type is as follows.
\begin{enumerate}
	\item Normalize $h$ such that it defines a probability density.
	\item Show that the sequences of probability measures $\nu_t$ defined by
\begin{equation}
	\nu_t(f)=\int_{\RR} f\left(\varphi\left(
	\begin{pmatrix} 1 & x \\ 0 & 1\end{pmatrix}
	\begin{pmatrix} \e^{-t/2} & 0 \\ 0 & \e^{t/2} \end{pmatrix}
	\right)\right) h(x) dx 
\end{equation}
	is tight. Then, by the Helly-Prokhorov theorem, it is relatively compact, i.e., every sequence of $\nu_t$ contains a convergent subsequence with weak limit $\nu$, say.
	\item Show that $\nu$ is invariant under a unipotent subgroup $U$; in the present case,
\begin{equation}
	U= \left\{\varphi\left(\begin{pmatrix} 1 & x \\ 0 & 1\end{pmatrix} \right)\right\}_{x\in\RR}.
\end{equation}
	\item Use a density argument to rule out measures concentrated on subvarieties (exploit the assumption that $\Gamma\backslash\Gamma H$ is dense in $\GamG$).
\end{enumerate}

As an application of Shah's theorem we give a proof of the statement in Remark \ref{rem2}, in fact a slightly more general version allowing for non-constant $h$. Recall that here $G=\SL(2,\RR)\times\RR^2$ and $\Gamma=\SL(2,\ZZ)\times\ZZ^2$.

\begin{cor}
Let $y\not\in\QQ$. For any bounded piecewise continuous $f:\GamG\to\RR$ and piecewise continuous $h:\RR\to\RR$ with compact support
\begin{equation}\label{soso}
	\lim_{t\to\infty} \int_{\RR} f(n_-(\alpha,y)\Phi^t) h(\alpha) d\alpha
	=\frac{1}{\mu(\GamG)} \int_{\GamG} f d\mu \int_\RR h(\alpha) d\alpha.
\end{equation}
\end{cor}

\begin{proof}
We define the embedding $\varphi:\SL(2,\RR)\to G$ by
\begin{equation}
	M \mapsto (1,(0,y)) (M,0) (1,(0,y))^{-1} .
\end{equation}
We need to show that 
\begin{equation}
	(\gamma,n)(1,(0,y)) (M,0) (1,(0,y))^{-1}
\end{equation}
is dense in $G$ as $\gamma,n,M$ vary over $\SL(2,\ZZ)$, $\ZZ^2$, $\SL(2,\RR)$, respectively. It is obviously sufficient to show this for
\begin{equation}
		(\gamma,n)(1,(0,y)) (M,0) = (\gamma M,(n_1,(y+n_2))M),
\end{equation}
and thus for $(M,(n_1,(y+n_2))\gamma^{-1} M)$.
It is however easy to see, using the irrationality of $y$, that $(n_1,(y+n_2))\gamma^{-1}$ is dense in $\RR^2$ (exercise). The completes the proof of the density.

Shah's theorem says now that
\begin{equation}
	\lim_{t\to\infty} \int_{\RR} \tilde f(n_-(\alpha,y)\Phi^t n_-(0,y)^{-1}) h(\alpha) d\alpha
	=\frac{1}{\mu(\GamG)} \int_{\GamG} \tilde f d\mu \int_\RR h(\alpha) d\alpha.
\end{equation}	
for all bounded, piecewise continuous $\tilde f$. Choosing the test function 
\begin{equation}
	\tilde f(g)=f(g n_-(0,y))
\end{equation}
which is left-$\Gamma$-invariant and bounded, piecewise continuous, if $f$ is (as assumed). This yields \eqref{soso}.
\end{proof}

\subsection{Equidistribution, unbounded test functions and diophantine conditions}

In some applications of Ratner's theorem, e.g., in questions of value distribution of quadratic forms \cite{Eskin98,Eskin05,pairI,pairII}, the test functions $f$ in the equidistribution theorems are no longer bounded. Under such circumstances the convergence of the integral can only be assured by assuming certain diophantine conditions. Without going into the intricate details for general $\GamG$, we will illustrate this phenomenon in the distribution of $m\alpha$ on $\TT$, which indeed may be viewed as a unipotent orbit on the homogeneous space $\ZZ\backslash\RR$. As mentioned earlier, it is well known that for $\alpha\notin\QQ$ the sequence is uniformly distributed mod 1. That is, for any bounded continuous function $f:\TT\to\RR$ we have
\begin{equation}
	\lim_{N\to\infty} \frac1N \sum_{m=1}^N f(m\alpha) = \int_{\TT} f(x) dx.
\end{equation}
Let us know formulate the analogous statement for test functions with a possible singularity at $x=0$.

It is convenient to identify $\TT$ with $[-1/2,1/2)$. For any $\beta\geq 0$ we define the class $K_\beta(\TT)$ of functions continuous on $\TT-\{0\}$, with the property that there is a constant $C>0$ such that
\begin{equation}
	|f(x)| \leq C |x|^{-\beta}, \qquad \text{for all $x\in[-1/2,1/2)$.}
\end{equation}

We say $\alpha\in\RR$ is {\em diophantine of type $\kappa$}
if there exists a constant $c>0$ such that
$$
\left| \alpha - \frac{p}{q} \right| > \frac{c}{q^\kappa}
$$
for all $p,q\in\ZZ$, $q>0$. The smallest possible value of $\kappa$ is $\kappa=2$ (achieved for quadratic surds, e.g., $\alpha=\sqrt2$), and it is well known that for any $\kappa>2$ there is a set of full Lebesgue measure of $\alpha$ that have type $\kappa$.

\begin{thm}
Let $\alpha$ be diophantine of type $\kappa$. Then, for any $f\in K_\beta(\TT)$ with $0\leq \beta<1/(\kappa-1)$,
\begin{equation}\label{mittach}
	\lim_{N\to\infty} \frac1N \sum_{m=1}^N f(m\alpha) = \int_{\TT} f(x) dx.
\end{equation}
\end{thm}

\begin{proof}
We split $f=f_+-f_-$ into positive and negative part, such that $f_\pm\geq 0$. Then $f_\pm\in K_\beta(\TT)$ and we may prove \eqref{mittach} for both $f_\pm$ separately. In the following we will thus assume that $f\geq 0$.

For any $\epsilon>0$ let us define
\begin{equation}
	f_\epsilon(x)=
	\begin{cases}
	f(x) & \text{if $|x|>\epsilon$}\\
	\min\{f(x), f(\epsilon)\} & \text{if $|x|\leq \epsilon$}
	\end{cases}
\end{equation}
and $g_\epsilon=f-f_\epsilon$. Note that $f_\epsilon\leq f$. By assumption,
\begin{equation}
	g_\epsilon(x)
		\begin{cases}
	\leq C |x|^{-\beta} & \text{if $|x|\leq \epsilon$}\\
	= 0 & \text{if $|x| \geq \epsilon$.}
	\end{cases}
\end{equation}
The function $f_\epsilon$ is bounded continuous, and hence by uniform distribution
\begin{equation} \label{feps}
	\lim_{N\to\infty} \frac1N \sum_{m=1}^N f_\epsilon(m\alpha) = \int_{\TT} f_\epsilon(x) dx
= \int_{\TT} f(x) dx - O(\epsilon^{1-\beta}).
\end{equation}
Since $f_\epsilon\leq f$, this implies the lower bound
\begin{equation}
	\liminf_{N\to\infty} \frac1N \sum_{m=1}^N f(m\alpha) 
\geq \int_{\TT} f(x) dx - O(\epsilon^{1-\beta}) . 
\end{equation}

As to the upper bound,
\begin{equation}
	\frac1N \sum_{m=1}^N g_\epsilon(m\alpha)
	\leq \frac{C}{N} \sum_{m=1}^N \frac{\chi_{(0,\epsilon]}(\|m\alpha\|)}{\|m\alpha \|^{\beta}}
\end{equation}
where $\|\,\cdot\,\|$ denotes the distance to the nearest integer.
Using the dyadic decomposition of the unit interval, we find
\begin{equation}
\begin{split}
	\frac{1}{N} \sum_{m=1}^N \frac{\chi_{(0,\epsilon]}(\|m\alpha\|)}{\|m\alpha \|^{\beta}}
	& = \frac{1}{N} \sum_{j=0}^\infty \sum_{m=1}^N \frac{\chi_{(\epsilon 2^{-(j+1)},\epsilon 2^{-j}]}(\|m\alpha\|)}{\|m\alpha \|^{\beta}} \\
	& < \frac{1}{N\epsilon^\beta} \sum_{j=0}^\infty 2^{\beta(j+1)} \sum_{m=1}^N \chi_{(\epsilon 2^{-(j+1)},\epsilon 2^{-j}]}(\|m\alpha\|) \\
	& \leq \frac{2B}{\epsilon^\beta} \sum_{j=0}^\infty 2^{\beta(j+1)} (\epsilon 2^{-(j+1)})^{\frac{1}{\kappa-1}} \quad\text{(for some $B>0$)}\\
	& =O(\epsilon^{\frac{1}{\kappa-1}-\beta}).
\end{split}
\end{equation}
The inequality before the last follows from Lemma \ref{lastLemma} below. 
We conclude that
\begin{equation}\label{gg}
	\frac1N \sum_{m=1}^N g_\epsilon(m\alpha) 
	= O(\epsilon^{\frac{1}{\kappa-1}-\beta})
\end{equation}
Therefore
\begin{equation}
\begin{split}
	\limsup_{N\to\infty} \frac1N \sum_{m=1}^N f(m\alpha)
	& =  \limsup_{N\to\infty} \frac1N \sum_{m=1}^N [f_\epsilon(m\alpha)+g_\epsilon(m\alpha)]\\
	& \leq \int_{\TT} f(x) dx + O(\epsilon^{1-\beta}) + O(\epsilon^{\frac{1}{\kappa-1}-\beta}),
\end{split}
\end{equation}
in view of \eqref{feps} and \eqref{gg}.

Since $\epsilon>0$ can be arbitrarily small, the limsup and liminf must coincide.
\end{proof}

The following lemma is used in the preceding proof.

\begin{lem}\label{lastLemma}
Let $\alpha$ be diophantine of type $\kappa$. Then there is a constant $B>0$ such that, for any interval $[x_0,x_0+\ell]$,
\begin{equation}\label{numb}
	\#\{ m=1,\ldots,N: m\alpha \in [x_0,x_0+\ell] +\ZZ \} 
	\leq 
	\begin{cases}
	0 & \text{if $N^{\kappa-1}\ell< c$} \\
	B N \ell^{1/(\kappa-1)} & \text{otherwise.} 
	\end{cases} 
\end{equation}
\end{lem}

\begin{proof}
Define $T=1/\ell$.
Let us divide the counting into blocks of the form
\begin{equation}
	\#\{ m_0 < m \leq m_0 + T^{1/(\kappa-1)} : m\alpha \in [x_0,x_0+\ell] +\ZZ \} ,
\end{equation}
The number of such blocks contributing to \eqref{numb} is less than $O(N T^{-1/(\kappa-1)}+1)$.

The gaps between elements of the sequence $m\alpha$ mod 1, $m_0< m\leq m_0+ T^{1/(\kappa-1)}$, are of the form $n\alpha$ mod 1, with $|n|<2T^{1/(\kappa-1)}$. By the diophantine condition, the gaps therefore have seize at least 
\begin{equation}
	\| n\alpha \| \geq \frac{c}{|n|^{\kappa-1}} > \frac{c}{2^{\kappa-1} T}. 
\end{equation}
An interval of size $\ell=1/T$ can hence at most contain a bounded number of elements. 
Hence
\begin{equation}
	\#\{ m_0 < m \leq m_0 + T^{1/(\kappa-1)} : m\alpha \in [x_0,x_0+\ell] +\ZZ \} \leq B'
\end{equation}
for some constant $B'>0$ independent of $m_0,x_0,\ell$. Recall that there were at most $N T^{-1/(\kappa-1)}+1$ such blocks, and this yields the upper bound in the second alternative. 

The first alternative is easily proven since the minimum gap size for the full sequence $m=1,\ldots,N$ is at least $c/(2N)^{\kappa-1}$. 
\end{proof}


\begin{thebibliography}{99}

\bibitem{Einsiedler06}
M. Einsiedler, A. Katok and E. Lindenstrauss,
Invariant measures and the set of exceptions to Littlewoods conjecture, to appear in {\em Annals of Math.}

\bibitem{Elkies04}
N.D. Elkies and C.T. McMullen, Gaps in ${\sqrt n}\bmod 1$ and ergodic theory.  {\em Duke Math. J.} {\bf 123}  (2004) 95--139.

\bibitem{Eskin98}
A.\,Eskin, G.\,Margulis and S.\,Mozes, Upper bounds and
asymptotics in a quantitative version of the Oppenheim conjecture,
{\em Ann. of Math.} {\bf 147} (1998) 93-141.

\bibitem{Eskin05}
A.\,Eskin, G.\,Margulis and S.\,Mozes, 
Quadratic forms of signature (2,2) 
and eigenvalue spacings on rectangular 2-tori, {\em Ann. of Math.} {\bf 161} (2005) 679-725.

\bibitem{Kleinbock99}
D. Kleinbock, Badly approximable systems of affine forms.  {\em J. Number Theory} {\bf 79}  (1999) 83-102. 

\bibitem{Lindenstrauss06}
E. Lindenstrauss,  Invariant measures and arithmetic quantum unique ergodicity,  {\em Ann. of Math.} {\bf 163}  (2006) 165-219.

\bibitem{npoint}
J.\,Marklof,
The $n$-point correlations between values of a linear form, 
with an appendix by Z.~Rudnick,
{\em Ergod. Th. Dyn. Sys.} {\bf 20} (2000) 1127-1172.

\bibitem{pairI}
J.\,Marklof, Pair correlation densities of inhomogeneous quadratic forms,
{\em Ann. of Math.} {\bf 158} (2003) 419-471.

\bibitem{pairII}
J.\,Marklof, Pair correlation densities of inhomogeneous quadratic forms II,
{\em Duke Math. J.} {\bf 115} (2002) 409-434; Correction,
{\em ibid.} {\bf 120} (2003) 227-228.

\bibitem{energy}
J. Marklof, Energy level statistics, lattice point problems 
and almost modular functions, in P. Cartier; B. Julia; P. Moussa; P. Vanhove (Editors):
{\em Frontiers in Number Theory, Physics and Geometry. Volume 1:
On random matrices, zeta functions and dynamical systems}, Springer, 2006, pp. 163-181.

\bibitem{Morris05}
D.W. Morris, {\em Ratner's theorems on unipotent flows.} Chicago Lectures in Mathematics. University of Chicago Press, Chicago, IL, 2005.

\bibitem{Rudnick98}
Z. Rudnick and P. Sarnak,
The pair correlation function of fractional parts of polynomials,
{\em Comm. Math. Phys.} \textbf{194} (1998) 61-70.

\bibitem{Rudnick02}
Z. Rudnick and A. Zaharescu, 
The distribution of spacings between fractional parts of lacunary sequences.
{\em Forum Math.} {\bf 14} (2002) 691-712.

\bibitem{Shah96}
N.A.\,Shah, Limit distributions of expanding 
translates of certain orbits on homogeneous spaces,
{\em Proc. Indian Acad. Sci., Math. Sci.} {\bf 106} (1996) 105-125.

\bibitem{Slater67}
N.B.\,Slater,
Gaps and steps for the sequence $n\theta \mod 1$, 
{\em Proc. Cambridge Philos. Soc.} {\bf 63} (1967) 1115-1123.

\bibitem{Strombersson05}
A. Str\"ombergsson and A. Venkatesh, Small solutions to linear congruences and Hecke equidistribution.  {\em Acta Arith.}  {\bf 118}  (2005) 41-78.

\end{thebibliography}
\end{document}